\documentclass[a4paper,11pt,reqno]{amsart}
\usepackage{amssymb,amsmath,latexsym,amsthm}

\newtheorem{theorem}{Theorem}[section]

\theoremstyle{definition}

\newtheorem*{remark}{Remark}
\newtheorem*{example}{Example}

\numberwithin{equation}{section}

\def\n{{\mathbb{N}}}

\def\z{{\mathbb{Z}}}
\def\ee{{\mathcal{E}}}
\def\cc{{\mathcal{C}}}
\def\pp{{\mathcal{P}}}
\def\uu{{\mathcal{U}}}
\def\ww{{\mathcal{W}}}

\begin{document}

\title[Absolute pseudoprimes with any number of prime factors]{New
polynomials producing absolute pseudoprimes with any number of prime
factors}

\author[Ken Nakamula]{Ken NAKAMULA}
\address{Ken Nakamula\\
Department of Mathematics and Information Sciences\\
Tokyo Metropolitan University\\
Minami-Osawa, Hachioji\\
192-0397 Tokyo, Japan}
\email{nakamula@tnt.math.metro-u.ac.jp}

\author[Hirofumi Tsumura]{Hirofumi TSUMURA}
\address{Hirofumi Tsumura\\
Department of Mathematics and Information Sciences\\
Tokyo Metropolitan University\\
Minami-Osawa, Hachioji\\
192-0397 Tokyo, Japan}
\email{tsumura@comp.metro-u.ac.jp}

\author[Hiroaki Komai]{Hiroaki KOMAI}
\address{Hiroaki Komai\\
Toa Gakuen High School\\
5-44-3 Kamitakada, Nakano-ku\\
164-0002 Tokyo, Japan}
\email{hiroaki-k.80@cup.ocn.ne.jp}

\thanks{
This research is partially supported by the JSPS Grand-in-Aid for
Scientific Research No.\ 16340011.
}

\maketitle

\begin{abstract}
In this paper, we introduce a certain method to construct polynomials
producing many absolute pseudoprimes.  By this method, we give new
polynomials producing absolute pseudoprimes with any fixed number of
prime factors which can be viewed as a generalization of Chernick's
result.  By the similar method, we give another type of polynomials
producing many absolute pseudoprimes.  As concrete examples, we
tabulate the counts of such numbers of our forms.
\end{abstract}

\bigskip

\section{Introduction} \label{S-1} 

A composite positive integer $N$ is called the \textit{absolute
pseudoprime} if $a^{N-1} \equiv 1$ (mod $N$) for any integer $a$ prime
to $N$. This is also called the \textit{Carmichael number}. The
following criterion is essential:

\smallskip

{\bf Korselt's criterion}. A composite odd positive integer $N$ is the
absolute pseudoprime if and only if $N$ is squarefree and $p-1$
divides $N-1$ for every prime $p$ dividing $N$.

\smallskip

From this, it easily follows that the number of prime divisors of the
absolute pseudoprime is at least three.

In \cite{AGP}, Alford, Granville and Pomerance showed that there are
infinitely many absolute pseudoprimes. Furthermore L{\H o}w and
Niebuhr \cite{LN} introduced a new algorithm for constructing absolute
pseudoprimes with a large number of prime factors. However, it is
still open that there are infinitely many $k$-component absolute
pseudoprimes for each $k \geq 3$. Hence we aim to produce many
$k$-component absolute pseudoprimes. Furthermore it seems to be
important to consider the polynomials producing many absolute
pseudoprimes with $k$ prime factors which is the product of distinct
$k$ linear polynomials such as
\[
U_{k}(m)=U_{k}\left( m;\{
\alpha_i\},\{\beta_i\}\right)=\prod_{i=1}^{k}(\alpha_i m+\beta_i)\ \ \
(m \in \n),
\]
where $\alpha_i \in \n$ and $\beta_i \in \z$ $(1 \leq i \leq k)$ which satisfy the congruences $U_k(m) \equiv 1$ (mod $\alpha_{i}m+\beta_i-1$) for any $m \in \mathbb{N}$ $(1 \leq i \leq k)$. This has been already considered by Chernick \cite{Ch}. Chernick called $U_{k}(m)$ the universal form. For example, he constructed 
\begin{equation}
U_k(m)=(6m+1)(12m+1)\prod_{i=1}^{k-2}\left( 9 \cdot 2^{i}m+1\right) \hspace{0.5in}(m \in \mathbb{N}) \label{ee-0}
\end{equation}
for $k \in \mathbb{N}$ with $k \geq 3$. Indeed, if $k \geq 4$ and $m=2^{k-4}M$ then $U_k(2^{k-4}M)$ is a absolute pseudoprime when all factors in the right-hand side of (\ref{ee-0}) are simultaneously prime numbers. He further gave a certain algorithm to construct a $k$-component universal form from a $k$-component absolute pseudoprime. However, by the Chernich method, it is impossible to construct many $k$-component universal forms because we cannot obtain even one $k$-component absolute pseudoprime for an arbitrary $k$. 

From the viewpoint of what is called the $k$-tuple prime conjecture (see \cite{CP} Chapter 1), it seem to be natural that $U_k(m)$ produces infinitely many $k$-component absolute pseudoprimes. Furthermore Granville and Pomerance \cite{GP} 
considered this deeply and gave the general theory about estimation of the number of $k$-component absolute pseudoprimes under the Hardy-Littlewood conjecture. As they mentioned in \cite{GP} Section 2, it is sufficient to consider the case $b_i=1$ for each $i$ in order to construct many $k$-component absolute pseudoprimes. Namely we consider the form
\begin{equation}
U_{k}(m)=\prod_{i=1}^{k}(\alpha_i m+1). \label{ee-1}
\end{equation}
We can see that (\ref{ee-1}) produces $k$-component absolute pseudoprimes if $U_k(m)$ satisfies the congruences $U_k(m) \equiv 1$ (mod $\alpha_{i}m$) for any $m \in \mathbb{N}$ $(1 \leq i \leq k)$. However they did not concretely 
construct $k$-component universal forms, namely polynomials producing many absolute pseudoprimes with $k$ prime factors. 

In the present paper, we give a certain sufficient condition that (\ref{ee-1}) is a polynomial producing many absolute pseudoprimes with $k$ prime factors 
(see Theorem \ref{Th-1}). By this result, we construct the polynomial $U_{k,l}(m)$ producing many absolute pseudoprimes which can be viewed as a generalization of Chernick's universal form (\ref{ee-0}). In particular when $l=3$, we see that $U_{k,3}(m)$ coincides with (\ref{ee-0}). Similarly, we give another type of polynomial $V_k(m)$ producing many absolute pseudoprimes with $k$ prime factors. As concrete examples, we tabulate the counts of such numbers by using the method similar to Dubner's one. Indeed, Dubner turned the method of Hardy and Littlewood precisely (see \cite{HL}), and tabulated the counts of absolute pseudoprimes of the form $\uu_{3,3}(M)$ (see \cite{Du}). We make use of his method, and tabulate the counts of absolute pseudoprimes of the form $\uu_{4,4}(M)$, $\uu_{5,5}(M)$ and $\ww_{4}(3M)$. 

\bigskip 

\section{Polynomials producing many absolute pseudoprimes} \label{S-2}

First we give the following theorem.

\smallskip

\begin{theorem}\label{Th-1} Let ${\bf a}=\{ a_1,a_2\ldots,a_{r}\} \subset \mathbb{N}$ with $a_1<\cdots<a_{r}$ and $r \geq 3$, which satisfy that
\begin{align}
& a_1+a_2+\cdots+a_{r-1}=a_{r}; \label{e-2} \\
& a_{j}\,|\,2a_{r}\hspace{0.5in}(1 \leq j \leq r); \label{e-3} \\
& {\rm GCD}(a_1,a_2,\ldots,a_r)=1. \label{e-4}
\end{align}
Put
\begin{equation}
U_{k}(m;{\bf a})=\prod_{\nu=1}^{r}(2a_{r}a_{\nu} m+1) \cdot \prod_{i=1}^{k-r}\left(2^{i+1}a_{r}^2m+1\right) \label{e-6}
\end{equation}
for $k \in \mathbb{N}$ with $k \geq r$. Suppose $m \in \mathbb{N}$
with $2^{k-r-1}\,|\,m$ when $k>r$ and $m \in \mathbb{N}$ is arbitrary
when $k=r$, and put $m=2^{k-r-1}M$ $(resp.\ m=M)$ when $k>r$ $(resp.\
k=r)$. Then $U_{k}(2^{k-r-1}M;{\bf a})$ $(resp.\ U_{k}(M;{\bf a}))$ is
the polynomial producing many absolute pseudoprimes with $k$ prime
factors when $k>r$ $(resp.\ k=r)$.
\end{theorem}

\begin{proof}
First we prove the case where $k=r$.  We can write $U_{r}(m;{\bf
a})=\sum_{\mu=0}^{r} C_{\mu}m^{\mu}$ as a polynomial in $m$ of degree
$r$.  We immediately check that $C_{0}=1$ and
\begin{equation}
C_{\mu}=(2a_{r})^{\mu}\sum_{1 \leq i_1<\cdots<i_{\mu} \leq r} a_{i_1}\cdots a_{i_{\mu}}\hspace{0.5in}(1 \leq \mu \leq r). \label{e-7}
\end{equation}
In particular, it follows from (\ref{e-2}) that $C_1=2a_{r}(a_1+\cdots+a_{r})=4a_{r}^2$. From (\ref{e-3}), we have $C_1 \equiv 0$ (mod $2a_{r}a_{j})$, namely
\begin{equation}
C_{1}m \equiv 0\ (\textrm{mod}\ 2a_{r}a_{j}m) \hspace{0.5in}(1 \leq j \leq r,\ m \in \mathbb{N}). \label{e-8}
\end{equation}
Suppose $\mu \geq 2$. Considering each case when $i_{\mu}<r$ and $i_{\mu}=r$, we obtain
\begin{align*}
C_{\mu}& =(2a_{r})^{\mu} \sum_{i_1<\cdots<i_{\mu} < r} a_{i_1}\cdots a_{i_{\mu}} +(2a_{r})^{\mu} a_{r} \sum_{i_1<\cdots<i_{\mu-1}} a_{i_1}\cdots a_{i_{\mu-1}} \\
      & =(2a_{r})^{\mu}\bigg\{ \sum_{i_1<\cdots<i_{\mu} < r} a_{i_1}\cdots a_{i_{\mu}} +a_{r} \sum_{i_1<\cdots<i_{\mu-1}} a_{i_1}\cdots a_{i_{\mu-1}}\bigg\}.
\end{align*}
Since $a_{j}\,|\,2a_{r}$ $(1 \leq j\leq r)$, we obtain $C_{\mu} \equiv 0$ (mod $2a_{r}a_{j})$, namely 
\begin{equation}
C_{\mu}m^{\mu} \equiv 0\ (\textrm{mod}\ 2a_{r}a_{j}m)\hspace{0.5in}(1 \leq j\leq r,\ m \in \mathbb{N}). \label{e-10}
\end{equation}
Combining (\ref{e-8}) and (\ref{e-10}), we have
\[
U_{r}(m;{\bf a})-1 \equiv \sum_{\mu=1}^{r} C_{\mu}m^{\mu} \equiv 0\
(\textrm{mod}\ 2a_{r}a_{j}m)\hspace{0.3in}(1 \leq j\leq r,\ m \in
\mathbb{N}).
\]
Therefore $U_{r}(m;{\bf a})$ is the polynomial producing many absolute pseudoprimes.

Secondly we assume $k \geq r+1$. As well as the above argument, we write $U_{k}(m;{\bf a})=\sum_{\mu=0}^{k}D_{\mu}m^{\mu}$. From (\ref{e-6}), we obtain $D_{0}=1$ and 
\begin{equation}
D_{\mu}=(2a_{r})^{\mu}\ \sum_{p=r-1-k+\mu}^{\mu} \ \sum_{\{ I_{p},G_{p} \}}\left\{ \prod_{i \in I_{p}}a_{i} \cdot \prod_{g \in G_{p}}\left( 2^{g}a_{r}\right)\right\}\ \ \ \ (1 \leq \mu \leq k), \label{e-11}
\end{equation}
where the sum $\sum_{\{ I_{p},G_{p}\}}$ is taken over all
$I_{p} \subset \{1,2,\ldots,r-1\}$ such that $\sharp I_{p}=p$ 
and all $G_{p} \subset \{ 0,1,\ldots,k-r\}$ such that $\sharp G_{p}=\mu-p$. 

In order to prove that $U_{k}\left( 2^{k-r-1}M;{\bf a}\right)$ is the polynomial producing many absolute pseudoprimes, we have only to prove that 
\begin{align}
& D_{\mu}2^{(k-r-1)\mu}M^{\mu} \equiv 0\ (\textrm{mod}\ 2^{k-r}a_{r}a_{j}M)\ \ \ (1 \leq j \leq r-1); \label{e-12}\\
& D_{\mu}2^{(k-r-1)\mu}M^{\mu} \equiv 0\ (\textrm{mod}\ 2^{2k-2r}a_{r}^2M) \label{e-13}
\end{align}
for $\mu=1,2,\ldots,k$ and $M \in \mathbb{N}$. 

When $\mu=1$, it follows from (\ref{e-2}) that
\[
D_{1} =(2a_{r})\left
( \sum_{\nu=1}^{r}a_{\nu}+\sum_{j=1}^{k-r}2^{j}a_{r}\right)
=2^{k-r}2a_{r}^2.
\]
Hence, from (\ref{e-3}), we have 
\begin{align*}
& D_{1}2^{(k-r-1)}M \equiv 0\ (\textrm{mod}\ 2^{k-r}a_{r}a_{j}M)\ \ \ (1 \leq j \leq r-1); \\
& D_{1}2^{(k-r-1)}M \equiv 0\ (\textrm{mod}\ 2^{2k-2r}a_{r}^{2}M)
\end{align*}
for $M \in \mathbb{N}$. Hence (\ref{e-12}) and (\ref{e-13}) hold for $\mu=1$. 

When $\mu \geq 2$, we have 
\begin{equation*}
D_{\mu}=\sum_{p=r-1-k+\mu}^{\mu} \ \sum_{\{ I_{p},G_{p} \}}\left\{
(2a_{r})^{p}\prod_{i \in I_{p}}N_{i} \cdot \prod_{g \in G_{p}}\left(
2^{g+1}a_{r}^2\right)\right\}.
\end{equation*}
Hence we have $D_{\mu} \equiv 0$ (mod $4a_{r}^2$). By (\ref{e-3}), we
see that (\ref{e-12}) and (\ref{e-13}) hold for $\mu\geq 2$. Thus
$U_{k}\left( 2^{k-r-1}M;{\bf a}\right)$ is the polynomial producing
many absolute pseudoprimes with $k$ prime factors if $k \geq r+1$.
\end{proof}

\begin{remark}
If $r-1$ strictly increasing natural numbers satisfy the condition
\[
a_{j}\,|\,2(a_1+a_2+\cdots+a_{r-1})\hspace{0.5in}(1 \leq j \leq r-1),
\]
then we can define $a_r$ by (\ref{e-2}) and can get ${\bf a}$ in the
theorem, factoring out the GCD of ${\bf a}$.
\end{remark}

\smallskip 

\begin{example} \label{ex-1}
We consider the cases $(r,k)=(4,5)$, $(5,6)$ and $(6,7)$. Then we can immediately give $r$-tuple sequences $(a_1,a_2,\ldots,a_r)$ ($r=4,5,6$) which satisfy (\ref{e-2})-(\ref{e-4}), for example, 
\begin{align*}
& (a_1,a_2,a_3,a_4)=(2,3,10,15),\ (2,12,28,42);\\
& (a_1,a_2,a_3,a_4,a_5)=(3,22,30,110,165),\ (6,14,15,70,105);\\
& (a_1,a_2,a_3,a_4,a_5,a_6)=(9,12,14,28,63,126),\ (15,20,21,70,84,210).
\end{align*}
Corresponding to these sequences, we can construct $U_{k}(m;{\bf a})$ ($k=5,6,7$) determined by (\ref{e-6}), such as 
\begin{align*}
& ( 60 m+1)( 90 m+1)( 300 m+1)( 450 m+1)( 900 m+1), \\
& ( 168 m+1)( 1008 m+1)( 2352 m+1)( 3528 m+1)( 7056 m+1); \\
& (990 m+1)(7260 m+1)(9900 m+1)(36300 m+1)\\
& \hspace{1in} \times (54450 m+1)(108900 m+1), \\
& (1260 m+1)(2940 m+1)(3150 m+1)(14700 m+1)\\
& \hspace{1in} \times (22050 m+1)(44100 m+1);\\
& (2268 m+1)(3024 m+1)(3528 m+1)(7056 m+1) \\
& \hspace{1in} \times (15876 m+1)(31752 m+1)(63504 m+1),\\
& (6300 m+1)(8400 m+1)(8820 m+1)(29400 m+1) \\
& \hspace{1in} \times (35280 m+1)(88200 m+1)(176400m+1).
\end{align*}
Furthermore we can construct certain classes of these forms as follows. 
Suppose $l \in \mathbb{N}$ with $l \geq 3$. Then we can apply Theorem \ref{Th-1} with $r=l$, $a_1=1$, $a_2=2^{l-2}$ and $a_j=2^{j-3}\left( 2^{l-2}+1\right)$, because the conditions (\ref{e-2})-(\ref{e-4}) hold. For $k \in \mathbb{N}$ with $k \geq l$, we define $\mathcal{U}_{k,l}(m)=U_{k}\left(m;\{a_1,\ldots,a_{l}\}\right)$. Then
\begin{align}
\mathcal{U}_{k,l}(m)= & \left( 2^{l-2} \left( 2^{l-2}+1\right)m+1\right) \left( 2^{2l-4} \left( 2^{l-2}+1\right)m+1\right) \label{e-14} \\
& \ \times \prod_{i=1}^{k-2} \left( 2^{l+i-3} \left( 2^{l-2}+1\right)^2 m+1\right). \nonumber
\end{align}
When $k=l$, then $\mathcal{U}_{k,k}(m)$ is the polynomial producing many absolute pseudoprimes with $k$ prime factors. When $k \geq l+1$, putting $m=2^{k-l-1}M$, we see that $\mathcal{U}_{k,l}\left( 2^{k-l-1}M\right)$ is the polynomial producing many absolute pseudoprimes. Note that if $l=3$ then $\mathcal{U}_{k,3}(m)$ coincides with the Chernick form (\ref{ee-0}). When $(k,l)=(4,4),\,(5,5)$, we have
\begin{align}
& \uu_{4,4}(m)=(20m+1)(80m+1)(100m+1)(200m+1), \label{e-15} \\
& \uu_{5,5}(m)=(72m+1)(576m+1)(648m+1) \label{e-15-2}\\
& \hspace{1.5in} \times (1296m+1)(2592m+1). \notag 
\end{align}
\end{example}

\smallskip 

By the same consideration as in the proof of Theorem \ref{Th-1}, we give another polynomial producing many absolute pseudoprimes as follows. However, this form can not be derived from Theorem \ref{Th-1} directly. 

\smallskip 

\begin{example} \label{ex-2} For $k \geq 3$, we define
\begin{equation}
\mathcal{W}_{k}(m)=(6m+1)\prod_{i=1}^{k-2}\left( 4\cdot 3^{i}m+1\right) \cdot \left( 2\cdot 3^{k-1}m+1\right). \label{e-16}
\end{equation}
When $3^{k-3}\,|\,m$, putting $m=3^{k-3}M$, we see that $\mathcal{W}_{k}\left( 3^{k-3}M\right)$ is the polynomial producing many absolute pseudoprimes with $k$ prime factors. In order to prove this fact, we have only to check that
\begin{align}
& \mathcal{W}_{k}\left( 3^{k-3}M\right) \equiv 1\ \ \left( \textrm{mod}\ 4\cdot 3^{2k-5}M\right); \label{e-17} \\
& \mathcal{W}_{k}\left( 3^{k-3}M\right) \equiv 1\ \ \left( \textrm{mod}\ 2\cdot 3^{2k-4}M\right), \label{e-18}
\end{align}
Indeed, if we write $\mathcal{W}_{k}(m)=\sum_{\mu=0}^{k}\ E_{\mu}m^{\mu}$ then we can see that $E_{0}=1$, $E_{1}=4\cdot 3^{k-1}$, and 
\[
E_{\mu}\left( 3^{k-3}M\right)^{\mu} \equiv 0\ \ \left( \textrm{mod}\
4\cdot 3^{2k-4}M\right)\ \ \ (\mu \geq 2).
\]
Hence we see that (\ref{e-17}) and (\ref{e-18}) hold and $\mathcal{W}_{k}\left( 3^{k-3}M\right)$ is the polynomial producing many absolute pseudoprimes with $k$ prime factors. When $k=3$, $W_{3}(m)$ coincides with Chernick's form $(6m+1)(12m+1)(18m+1)$. When $k=4$, we have the polynomial producing many absolute pseudoprimes with $4$ prime factors
\begin{equation}
\mathcal{W}_{4}(3M)=(18M+1)(36M+1)(108M+1)(162M+1). \label{e-19}
\end{equation}
\end{example}

\ 

\section{Numerical results}
In \cite{Du}, Dubner turned the method of Hardy and Littlewood precisely (see \cite{HL}) and constructed a function for estimating the count of absolute pseudoprimes of the form $\uu_{3,3}(m)=(6m+1)(12m+1)(18m+1)$. We apply his method with $\uu_{4,4}$, $\uu_{5,5}$ and $\ww_4$. 

We recall Dubner's method to estimate the number of absolute pseudoprimes (see \cite{Du}\ \S 3). 
Denote by $\pp(N)$ the probability of $N$ being prime for $N \in \n$. By the Prime Number Theorem, we have 
\begin{equation}
\pp(N) \sim  \frac{1}{\log N}\ \ \ \ (N \to \infty). \label{e-20}
\end{equation}
On the other hand, it follows from the Mertens theorem  (see \cite{HW} $\S\, 22.8$) that
\begin{equation}
\prod_{p:\textrm{prime} \atop p \leq \sqrt{N}}\,\frac{p-1}{p} \sim \frac{2e^{-\gamma}}{\log N}\ \ \ (N \to \infty). \label{e-21}
\end{equation}
First, we consider $\uu_{4,4}(m)$. Let $u=q\cdot r\cdot s\cdot t$, where $q=20m+1$, $r=80m+1$, $s=100m+1$ and $t=200m+1$. By the Prime Number Theorem, the probability of $q$ being prime becomes 
\begin{equation}
P_{q}=\frac{2}{2-1}\cdot\frac{5}{5-1}\cdot\frac{1}{\log(20m+1)}=\frac{2.5}{\log(20m+1)},
\end{equation}
because $q$ is not divisible by $2$ or $5$. However, the primality of $r$ is affected if $q$ is prime, because $q=20m+1$ and $r=80m+1$, namely $r=4q-3$. Let $p$ is a prime with $7 \leq p \leq \sqrt{r}$. Then the condition that $q$ is prime shows that $r \not\equiv 4p-3 \equiv -3$ (mod $p$). Under this condition, if $r$ is prime then $r \not\equiv 0,p-3$ (mod $p$). This means that we need to consider the correction factor $\cc_{r}(p)$ defined by 
\[
\cc_{r}(p)=\frac{p}{p-1}\cdot
\frac{p-2}{p-1}=\frac{p(p-2)}{(p-1)(p-1)}.
\]
When $p=3$, the condition that $q$ is prime shows that $r \not\equiv 0$ (mod $3$). Hence we let
\begin{equation*}
\cc_{r}=\frac{2}{1}\cdot \frac{3}{2}\cdot \frac{5}{4}\cdot \prod_{p \geq 7}\frac{p(p-2)}{(p-1)(p-1)} \doteqdot 3.520865,
\end{equation*}
Then the probability of $r$ being prime becomes 
\begin{equation}
P_{r}=\cc_{r}\cdot \frac{1}{\log(80m+1)}\doteqdot \frac{3.520865}{\log(80m+1)} . \label{e-22}
\end{equation}
We see that $4s=5r-1$. For a prime $p$ with $7 \leq p \leq \sqrt{r}$, if $r \not\equiv 0,-3$ (mod $p$) then $s \not\equiv -4,-4^{-1}$ (mod $p$), where $4^{-1}$ is a inverse element of $4$ mod $p$. Note that $4 \not\equiv 4^{-1}$ (mod $p$), since $p>5$. Under the condition that both $q$ and $r$ are prime, if $s$ is prime then $s \not\equiv 0,-4,-4^{-1}$ (mod $p$). This means that we need to consider the correction factor $\cc_{s}(p)$ defined by 
\[
\cc_{s}(p)=\frac{p}{p-1}\cdot
\frac{p-3}{p-2}=\frac{p(p-3)}{(p-1)(p-2)}.
\]
When $p=3$, the condition that $r \not\equiv 0$ (mod $3$) means that $s \not\equiv -1$ (mod $3$) because $4s=5r-1$. Under the condition, $s$ is prime means $s \not\equiv 0,-1$ (mod $3$). So we let
\begin{equation*}
\cc_{s}=\frac{2}{1}\cdot \left(\frac{3}{2}\cdot\frac{1}{2}\right)\cdot \frac{5}{4}\cdot \prod_{p \geq 7}\frac{p(p-3)}{(p-1)(p-2)} \doteqdot 1.623609,
\end{equation*}
Then the probability of $s$ being prime becomes 
\begin{equation}
P_{s}=\cc_{s}\cdot \frac{1}{\log(100m+1)}\doteqdot \frac{1.623609}{\log(100m+1)} . \label{e-23}
\end{equation}
Furthermore we see that $t=2s-1$. For a prime $p$ with $7 \leq p \leq \sqrt{r}$, if $s \not\equiv 0,-4,-4^{-1}$ (mod $p$) then $t \not\equiv -1, -9,-2\cdot4^{-1}-1$ (mod $p$). Note that $1,9 \not\equiv 2\cdot4^{-1}+1$ (mod $p$), since $p>5$. Under the condition that both $q,r$ and $s$ are prime, if $t$ is prime then $t \not\equiv 0,-1, -9,-2\cdot4^{-1}-1$ (mod $p$). This means that we need to consider the correction factor $\cc_{s}(p)$ defined by 
\[
\cc_{s}(p)=\frac{p}{p-1}\cdot
\frac{p-4}{p-3}=\frac{p(p-4)}{(p-1)(p-3)}.
\]
When $p=3$, the condition $s \not\equiv 0,-1$ (mod $3$) means $s \not\equiv 0,-1$ (mod $3$) because $t=2s-1$. Under the condition, $t$ is prime means $t \not\equiv 0,-1$ (mod $3$). Hence we let
\begin{equation*}
\cc_{t}=\frac{2}{1}\cdot \frac{3}{2} \cdot \frac{5}{4}\cdot \prod_{p \geq 7}\frac{p(p-4)}{(p-1)(p-3)} \doteqdot 2.904708,
\end{equation*}
Then the probability of $s$ being prime becomes 
\begin{equation}
P_{t}=\cc_{t}\cdot \frac{1}{\log(200m+1)}\doteqdot \frac{2.904708}{\log(200m+1)} . \label{e-24}
\end{equation}
Hence the probability of $q,r,s$ and $t$ being prime simultaneously becomes
\begin{align}
P_{qrst}& =P_{q}P_{r}P_{s}P_{t} \label{e-25}\\
& \doteqdot \frac{41.511967}{\log(20m+1)\log(80m+1)\log(100m+1)\log(200m+1)}. \notag
\end{align}
Following Dubner's method introduced in \cite{Du} Section 3, we consider 
\begin{align*}
& \ee_1 (M)=41.511967\sum_{m=1}^{M}\frac{1}{\log(20m+1)\log(80m+1)} \\
& \hspace{1.5in} \ \times \frac{1}{\log(100m+1)\log(200m+1)},
\end{align*}
which gives an estimate for the number of such absolute pseudoprimes with $m \leq M$ for a given $M$. 
Define $a_m$ by
\[
\log(20m+1)\log(80m+1)\log(100m+1)\log(200m+1)=(\log(a_m\cdot m))^4.
\]
Then the estimate becomes
\[
\ee_1 (M) \sim 41.511967\int_{1}^{M}\frac{dm}{(\log(a_{_M}\cdot
m))^4}.
\]
Integrating by parts third times gives
\begin{align}
\ee_1 (M) \sim \frac{41.511967}{6a_{_M}} & \bigg\{ \int_{a_{_M}}^{a_{_M}M}\frac{dt}{\log t}-\frac{a_{_M}M}{\log(a_{_M}M)} \label{e-26} \\
&\ \ -\frac{a_{_M}M}{(\log(a_{_M}M))^2}-\frac{2a_{_M}M}{(\log(a_{_M}M))^3}\bigg\}. \notag
\end{align}
Note that the first term in the right-hand side of (\ref{e-26}) can be calculated by using the well-known logarithmic integral function $L_i(x)$. By (\ref{e-26}), we obtain the following table of theoretical count of $\uu_{4,4}(m)$. Note that $N_1(M)$ is the actual number of such absolute pseudoprimes with $m \leq M$. 

\ 

\begin{center}
\begin{tabular}{|r|r|r|r|} \hline
\multicolumn{1}{|c|}{$M$}  & $\ee_1(M)$ & $N_1(M)$ & $\ee_1(M)/N_1(M)$ \\ \hline
 & & & \\
\textrm{$10^3$} & 2 & 2 & 1.00000\\
\textrm{$10^4$} & 16 & 17 & 0.94118\\
\textrm{$10^5$} & 90 & 87 & 1.03448\\
\textrm{$10^6$} & 506 & 487 & 1.03901\\
\textrm{$10^7$} & 3021 & 2959 & 1.02095\\
\textrm{$10^8$} & 19143 & 18960 & 1.00965\\
\textrm{$10^9$} & 127204 & 126997 & 1.00163\\
 & & & \\ \hline
\end{tabular}
\end{center}
\begin{center}
Table 1: Count of $\uu_{4,4}(m)$
\end{center}

\ 

Using the similar consideration, we can give an estimate for the number of absolute pseudoprimes of the form 
\[
\uu_{5,5}(m)=(72m+1)(576m+1)(648m+1)(1296m+1)(2592m+1),
\]
and of 
\[
\mathcal{W}_{4}(3m)=(18m+1)(36m+1)(108m+1)(162m+1).
\]
Corresponding to these forms, we let
\begin{align*}
& \ee_2 (M)=263.428500 \sum_{m=1}^{M}\frac{1}{\log(72m+1)\log(576m+1)\log(648m+1)} \\
& \hspace{1.5in} \times \frac{1}{\log(1296m+1)\log(2592m+1)}, \\
& \ee_3 (M)=66.419105\sum_{m=1}^{M}\frac{1}{\log(18m+1)\log(36m+1)} \\
& \hspace{1.5in} \times \frac{1}{\log(108m+1)\log(162m+1)}.
\end{align*}
Then we obtain the following table of theoretical count of $\uu_{5,5}(m)$ and $\mathcal{W}_{4}(3m)$, where $N_2(M)$ and $N_3(M)$ are the actual numbers of such absolute pseudoprimes with $m \leq M$. 

\ 

\begin{center}
\begin{tabular}{|r|r|r|r|} \hline
\multicolumn{1}{|c|}{$X$}  & $\ee_2(M)$ & $N_2(M)$ & $\ee_2(M)/N_2(M)$ \\ \hline
 & & & \\
\textrm{$10^3$} & 1 & 2 & 0.50000\\
\textrm{$10^4$} & 2 & 5 & 0.40000\\
\textrm{$10^5$} & 19 & 22 & 0.86364\\
\textrm{$10^6$} & 105 & 107 & 0.98131\\
\textrm{$10^7$} & 596 & 616 & 0.96753\\
\textrm{$10^8$} & 3555 & 3516 & 1.01109\\
\textrm{$10^9$} & 22261 & 22163 & 1.00442\\
 & & & \\ \hline
\end{tabular}
\end{center}
\begin{center}
Table 2: Count of $\uu_{5,5}(m)$
\end{center}

\ 

\begin{center}
\begin{tabular}{|r|r|r|r|} \hline
\multicolumn{1}{|c|}{$X$}  & $\ee_3(M)$ & $N_3(M)$ & $\ee_3(M)/N_3(M)$ \\ \hline
 & & & \\
\textrm{$10^3$} & 7 & 10 & 0.70000\\
\textrm{$10^4$} & 30 & 33 & 0.90909\\
\textrm{$10^5$} & 155 & 149 & 1.04027\\
\textrm{$10^6$} & 862 & 824 & 1.04612\\
\textrm{$10^7$} & 5108 & 5116 & 0.99843\\
\textrm{$10^8$} & 32170 & 32077 & 1.00290\\
\textrm{$10^9$} & 212716 & 213075 & 0.99832\\
 & & & \\ \hline
\end{tabular}
\end{center}
\begin{center}
Table 3: Count of $\ww_{4}(3m)$
\end{center}


\begin{thebibliography}{100}

\bibitem {AGP} \textsc{W. R. Alford, A. Granville and C. Pomerance}, \textit{There are infinitely many Carmichael numbers}, {Ann. Math.} \textbf{140} (1994), 703-722.

\bibitem{Ch} \textsc{J. Chernick}, \textit{On Fermat's simple theorem}, {Bull. Amer. Math. Soc.} {\bf 45} (1939), 269-274.

\bibitem{CP} \textsc{R. Crandall and C. Pomerance}, Prime numbers, Springer-Verlag, New-York Berlin Heidelberg, 2001.

\bibitem{Du} \textsc{H. Dubner}, \textit{Carmichael numbers of the form $(6m+1)(12m+1)(18m+1)$}, {J. Integer Seq.} {\bf 5} (2002), Article 02.2.1, 1-8.

\bibitem{GP} \textsc{A. Granville and C. Pomerance}, \textit{Two contradictory conjectures concerning Carmichael numbers}, {Math. Comp.} {\bf 71} (2002), 883-908.

\bibitem{HL} \textsc{G. H. Hardy and J. E. Littlewood}, \textit{Some problems on partitio numerorum III}, On the expression of a number as a sum of primes, {Acta Math.} {\bf 44} (1923), 1-70.

\bibitem{HW} \textsc{G. H. Hardy and E. M. Wright}, An introduction to the theory of numbers, fifth edition, Oxford Univ. Press, 1979.

\bibitem{LN} \textsc{G. L{\H o}h and W. Niebuhr}, \textit{A new algorithm for constructing large Carmichael numbers}, {Math. Comp.} {\bf 65} (1996), 823-836.

\end{thebibliography}
\end{document}